\documentclass[12pt,a4paper,reqno]{amsart}
\usepackage{amsmath}
\usepackage{amsthm}
\usepackage{amsfonts}
\usepackage{amssymb}
\numberwithin{equation}{section}

\theoremstyle{plain}

\newtheorem{theorem}[subsection]{Theorem}
\newtheorem{proposition}[subsection]{Proposition}
\newtheorem{lemma}[subsection]{Lemma}
\newtheorem{corollary}[subsection]{Corollary}

\newtheorem{definition}[subsection]{Definition}

\theoremstyle{definition}

\renewcommand{\leq}{\leqslant}
\renewcommand{\geq}{\geqslant}

\newsavebox{\proofbox}
\savebox{\proofbox}{\begin{picture}(7,7)%
 \put(0,0){\framebox(7,7){}}\end{picture}}

\def\Z{\mathbb{Z}}
\def\R{\mathbb{R}}

\def\N{\mathbb{N}}

\begin{document}

\title[Maximal number of 3APs]{On the maximal number of three-term arithmetic progressions in subsets of $\Z/p\Z$}

\author{Ben Green}
\address{Centre for Mathematical Sciences\\
     Wilberforce Road\\
     Cambridge CB3 0WA\\
     England
}
\email{b.j.green@dpmms.cam.ac.uk}

\author{Olof Sisask}
\address{Department of Mathematics\\
     University of Bristol\\
     Bristol BS8 1TW\\
     England
}
\email{O.Sisask@dpmms.cam.ac.uk}

\thanks{The first author was a Clay Research Fellow while this work was carried out and gratefully acknowledges the support of the Clay Institute. The second author is funded by an EPSRC DTG through the University of Bristol. While this work was being carried out, he spent time at MIT and the University of Cambridge, and would like to thank both institutions for their kind hospitality.}


\begin{abstract}
Let $\alpha \in [0,1]$ be a real number. Ernie Croot \cite{croot1} showed that the quantity
\[ \max_{\substack{A \subseteq \Z/p\Z \\ |A| = \lfloor \alpha p\rfloor}} \frac{\#(\mbox{3-term arithmetic progressions in $A$})}{p^2}\] tends to a limit as $p \rightarrow \infty$ though primes. Writing $c(\alpha)$ for this limit, we show that
\[ c(\alpha) = \alpha^2/2\] provided that $\alpha$ is smaller than some absolute constant. In fact we prove rather more, establishing a structure theorem for sets having the maximal number of 3-term progressions amongst all subsets of $\Z/p\Z$ of cardinality $m$, provided that $m < cp$.
\end{abstract}

\maketitle

     \section{Introduction}
     
There are many papers in the additive combinatorics literature in which a study is made of arithmetic progressions inside rather arbitrary sets of integers or residues. Perhaps most famous amongst these is Roth's 1953 paper \cite{roth} in which it is established that if $\alpha > 0$, and if $A \subseteq \{1,\dots,N\}$ is a set with cardinality at least $\alpha N$, then $A$ contains a 3-term arithmetic progression (3AP) provided that $N > N_0(\alpha)$ is sufficiently large. A subsequent argument by Varnavides \cite{varnavides} deduces from Roth's theorem that there are in fact at least $f_3(\alpha) N^2$ such 3APs, for some $f_3(\alpha) > 0$.    

These results establish lower bounds on the \emph{minimum} number of 3APs inside a set. Our focus in this paper will be on the \emph{maximum} number of 3APs a set may have. 

Given a set $A$ inside some abelian group $G$ with no $2$- or $3$-torsion we write $T_3(A)$ for the number of $x,d \in G$ for which $x, x+d, x+2d \in A$. This count of three-term progressions is the most natural one in many ways. Note however that it does count ``trivial'' progressions for which $d = 0$ (though in our settings these will never make a substantial contribution). More importantly it counts each ``combinatorial'' progression twice -- for example $(5,7,9)$ is counted by $(x,d) = (5,2)$ and by $(x,d) = (9,-2)$. 

Our first result is for sets of integers. We determine $M_3(n)$, the maximum number of 3APs in a set of integers of size $n$, and we also describe the structure of sets which have the maximal number of 3APs.
\begin{definition}[Extremal sets]\label{def2.1}
Suppose that $k,m \in \N \cup \{0\}$. Then we write
\[ E(k,m) := \{-k-2m,\dots,-k-2,-k,-k+1,\dots,-1,0,1,\dots,k-1,k,k+2,\dots,k+2m\}\]
and
\[ F(k,m) := \{-k-2m+2,\dots,-k-2,-k,-k+1,\dots,-1,0,1,\dots,k-1,k,k+2,\dots,k+2m\}.\]
\end{definition}
Note that each such set can be viewed as a disjoint union of two arithmetic progressions, each of common difference 2, and that 
\[ |E(k,m)| = 2k + 2m + 1 \qquad \mbox{and} \qquad |F(k,m)| = 2k + 2m.\]
\begin{theorem}\label{max3APInt}
We have $M_3(n) = \lceil n^2/2\rceil$. Furthermore if $|A| = n$ and $T_3(A) = M_3(n)$ then $A$ is an affine image of some set $E(k,m)$ or $F(k,m)$.\end{theorem}
It is interesting, and perhaps more natural, to ask questions about arithmetic progressions for subsets of \emph{finite} groups such as $\Z/N\Z$. One reason for this is the following simple result, which has no analogue when working in $\Z$.

\begin{lemma}[Progressions in a set and its complement]
Suppose that $G$ is a group with no $2$- or $3$-torsion, and that $A \subseteq G$ has cardinality $\alpha|G|$. Then
\[ T_3(A) + T_3(A^c) = (1 - 3\alpha + 3\alpha^2)|G|^2.\]
\end{lemma}     
\proof If $f_1,f_2,f_3 : G \rightarrow \R$ are any three functions then we write (by a slight abuse of notation)
\[ T_3(f_1,f_2,f_3) := \sum_{x,d \in G} f_1(x) f_2(x+d)f_3(x+2d).\]
Note that $T_3$ is a trilinear form and that 
\[ T_3(A) = T_3(1_A,1_A,1_A)\] for any set $A$.
Now we have
\[ T_3 (1_{A^c}, 1_{A^c}, 1_{A^c}) = T_3(1 - 1_A, 1 - 1_A, 1 - 1_A),\]
which may obviously be expanded as a sum of eight terms. It is clear that any term with precisely one $1_A$ is equal to $\alpha |G|^2$, and any with two copies of $1_A$ is equal to $\alpha^2 |G|^2$. The result follows quickly.\endproof

Thus the maximal number of 3APs in a subset of $\Z/N\Z$ of size $n$ is equal to the minimal number of 3APs in a subset of size $N - n$. Croot \cite{croot1} studied these problems and proved the following pleasant result.

\begin{proposition}[Croot]\label{crootprop} Let $1 \leq n \leq N$ and write \[ M_3(n,N) := \max_{\substack{A \subseteq \Z/N\Z \\ |A| = n}} T_3(A)\qquad \mbox{and} \qquad m_3(n,N) := \min_{\substack{A \subseteq \Z/N\Z \\ |A| = n}} T_3(A).\]
 Suppose that $\alpha \in [0,1]$ is a fixed real number. Then
\[ M_3(\alpha) := \lim_{\substack{N \rightarrow \infty \\ N \mbox{\scriptsize \textup{prime}}}} \frac{M_3(\lfloor \alpha N\rfloor, N)}{N^2}\qquad \mbox{and}
\qquad m_3(\alpha) := \lim_{\substack{N \rightarrow \infty \\ N \mbox{\scriptsize \textup{prime}}}} \frac{m_3(\lfloor \alpha N\rfloor, N)}{N^2}\]
exist. Furthermore we have the relation
\[ m_3(\alpha) + M_3(1-\alpha) = 1 - 3\alpha + 3\alpha^2.\]
\end{proposition}

Our main result is the following, which relates $M_3(n,N)$ to $M_3(n)$ provided that $N$ is sufficiently large in terms of $n$.
\begin{theorem}[Maximal number of 3APs modulo a prime]\label{mainthm}
There is an absolute constant $c > 0$ with the following property. If $n$ and $N$ are integers with $N$ prime and $n \leq cN$ then $M_3(n,N) = M_3(n) = \lceil n^2/2\rceil$. Furthermore the only subsets of $\Z/N\Z$ with cardinality $n$ and the maximal number of 3APs are affine images of the sets of Definition \ref{def2.1}.
\end{theorem}
This immediately implies the following result concerning Croot's function $M_3(\alpha)$.
\begin{corollary}
Suppose that $\alpha < c$. Then $M_3(\alpha) = \alpha^2/2$.
\end{corollary}
We make some miscellaneous further observations on Croot's function in \S \ref{miscell-sec}.

\section{Arithmetic progressions in sets of integers}

Our objective in this section is to give the (straightforward) proof of Theorem \ref{max3APInt}, which gave an evaluation of $M_3(n)$, the maximal number of 3APs that a set of $n$ integers may contain. It also classified those sets with the maximal number of 3APs as being affine images of one of the special types of set $E(k,m), F(k,m)$ defined in Definition \ref{def2.1}.

\emph{Proof of Theorem \ref{max3APInt}.} Let $A \subseteq \Z$ be a set with $|A| = n$ and $T_3(A) = M_3(n)$. Write $a_1, \ldots, a_n$ for the elements of $A$, listed in increasing order. Given an index $j$, the element $a_j$ can occur as the mid-point of at most $\min(j-1,n-j)$ increasing 3APs (that is, 3APs $(x,x+d,x+2d)$ with $d > 0$). Counting each such progression twice (for it may also be realised as $(x+2d, (x+2d) - d, (x+2d) - 2d)$) and remembering to include the trivial progressions $(x,x,x)$ we obtain 
\[ T_3(A) \leq n + 2\sum_{j=1}^n \min(j-1,n-j) = \lceil n^2/2 \rceil.\]
Equality holds if and only if every point $a_j$ is the mid-point of exactly $\min(j-1, n-j)$ increasing 3APs in $A$, and a short check confirms that this is indeed the case when $A$ belongs to one of the two families $E(k,m)$ and $F(k,m)$. 

It is only a little harder to show that these are, up to affine equivalence, the \emph{only} examples where equality holds. 

\emph{Case 1: $n$ is odd.} Write $n = 2t + 1$
Now $a_{t+1}$ must be the midpoint of $t$ increasing 3APs, which must therefore be $(a_i,a_{t+1},a_{2t+2-i})$ for $i = 1,\dots,t$. Consider now the point $a_{t}$, which must be the mid-point of exactly $t-1$ increasing 3APs. Noting in view of the preceding that $(a_i,a_t,a_{2t+1})$ is not a progression, we see that these progressions must be precisely \[ (a_{t-i},a_t,a_{t+i}) \qquad 1 \leq i \leq k-1;\]
\[ (a_{t-i},a_t,a_{t+i+1}) \qquad k \leq i \leq t-1\] for some $k$, $1 \leq k \leq t$. It is easy to check that this forces $A$ to be an affine image of $E(k,t-k)$.

\emph{Case 2: $n$ is even.} Write $n = 2t$. Clearly it is not possible for both of the triples $(a_1,a_t,a_{2t})$ and $(a_1,a_{t+1},a_{2t})$ to be 3APs. By sending $A$ to $-A$ if necessary we may assume that $a_1,a_t$ and $a_{2t}$ do not lie in arithmetic progression. Now the point $a_t$ must be the midpoint of $t-1$ progressions, which must therefore be $(a_i,a_t,a_{2t - i})$ for $i = 1,\dots,t-1$. The point $a_{t+1}$ must also lie in $t-1$ arithmetic progressions, which must be precisely
\[ (a_{t+1 - i}, a_{t+1}, a_{t+1+i}) \qquad 1 \leq i \leq k-1;\]
\[ (a_{t-i}, a_{t+1}, a_{t+1+i}) \qquad k \leq i \leq t-1\] for some $k$, $1 \leq k \leq t$. 
One may check that these conditions force $A$ to be an affine image of $F(k,t-k)$.\endproof

\emph{Remark.}
This proof intrinsically uses the fact that $\Z$ is an ordered group, and so fails in $\Z/N\Z$. One may also prove the result by induction, using the fact that either the smallest or the largest element of $A$ cannot be involved in too many 3APs; again, this uses the ordering of the integers in an essential way.

\section{A rough structure theorem for arbitrary additive sets}

In this section, and for the rest of the paper, the letters $C$ and $c$ will denote positive absolute constants which may vary from line to line.

A key ingredient of our work is Proposition \ref{structure-theorem} below, in which an arbitrary additive set $A$ is decomposed into $k$ disjoint ``additively structured'' parts $A_1,\dots,A_k$ plus a leftover set $A_0$, in such a way that there is little ``additive communication'' between different sets $A_i,A_j$. Our result is very close in spirit to a result of Elekes and Ruzsa \cite{elekes-ruzsa}, but does not seem to follow directly from it. 

Before stating the result, we recall the definition and basic properties of sumsets and \emph{additive energy}. For more details, \cite[Chapter 2]{tao-vu} may be consulted.

If $A,B$ are subsets of an abelian group then we write $A - B := \{a - b: a \in A, b \in B\}$.
If $\lambda \in \N$ then we write $\lambda \cdot A := \{\lambda a : a \in A\}$ and $\lambda A := \{a_1 + \dots + a_{\lambda} : a_1,\dots, a_{\lambda} \in A\}$. We define the \emph{additive energy} between $A$ and $B$ to be the quantity
\[ E(A,B) := \# \{(a_1,b_1,a_2,b_2) \in A \times B \times A \times B: a_1 + b_1 = a_2 + b_2\}. \] 
Write $\delta[A] := |A - A|/|A|$ for the growth of $A$ under the differencing operation.

\begin{lemma}[Basic properties of the additive energy]\label{basic-add-energy} Suppose that $A,B$ are two sets in an abelian group $G$. Let $\eta \in (0,1]$ be a real parameter. 
\begin{enumerate}
\item $E(A,B)$ is bounded by all three of the quantities $|A|^2|B|,|B|^2|A|$ and $|A|^{3/2}|B|^{3/2}$.
\item There is some $x$ such that $|A \cap (B + x)| \geq E(A,B)/|A||B|$.
\item $E(A,B) \geq |A|^2|B|^2/|A\pm B|$.
\item Suppose that $A,B$ are two additive sets with $\delta[A] = K_A$, $\delta[B] = K_B$ and $E(A,B) \geq \eta|A|^{3/2}|B|^{3/2}$. Then $\delta[A \cup B] \leq 4K_AK_B/\eta$.
\item \textup{(Balog-Szemer\'edi-Gowers theorem)} Suppose that $A,B$ are two additive sets with $E(A,B) \geq \eta |A|^{3/2}|B|^{3/2}$. Then there are sets $A' \subseteq A, B' \subseteq B$ such that $|A'| \geq c\eta^C|A|$, $|B'| \geq c\eta^C|B|$ and $|A'-B'| \leq C\eta^{-C}|A'|^{1/2}|B'|^{1/2}$.
\item Suppose that $G$ has no elements of order $\leq L$, and suppose that $0 < \lambda \leq L$. Suppose that $E(A,B) \leq \eta |A|^{3/2}|B|^{3/2}$. Then $E(\lambda \cdot A, B) \leq (C\eta)^{c/L}|A|^{3/2}|B|^{3/2}$.
\end{enumerate}
\end{lemma}
\proof (i) The first two bounds are immediate, and the third follows from the first two.

(ii) follows immediately from the chain of inequalities
\[ E(A,B) = \sum_x |A \cap (B + x)|^2 \leq \sup_x |A \cap (B + x)|\sum_x |A \cap (B + x)| = |A||B|\sup_x |A \cap (B+x)|.\]
(iii) follows from the Cauchy-Schwarz inequality. Writing $r(x)$ for the number of representations of $x$ as $a + b$, we have
\[ E(A,B) = \sum_x r(x)^2 \geq \frac{1}{|A+B|} \big( \sum_x r(x) \big)^2 = \frac{|A|^2|B|^2}{|A+B|}.\]
An essentially identical argument works for $A - B$.

(iv) Using part (ii), choose $x$ such that $S = A \cap (B + x)$ has size at least $\eta |A|^{1/2}|B|^{1/2}$. Since $S \subseteq A$ we have
\[ |A-S| \leq |A-A| \leq K_A|A|,\]
and since $S \subseteq B +x$ we have
\[ |B-S| \leq |B - B + x| = |B-B| \leq K_B|B|.\]
Thus by an instance of the Ruzsa triangle inequality (cf. \cite[Lemma 2.6]{tao-vu}) we have
\[ |A-B| \leq \frac{|A-S||B-S|}{|S|} \leq \frac{K_AK_B}{\eta}|A|^{1/2}|B|^{1/2}.\]
It follows that 
\[ |(A\cup B) - (A \cup B)| \leq K_A |A| + K_B |B| + \frac{2K_AK_B}{\eta}|A|^{1/2}|B|^{1/2}.\]
It is immediate from this that
\[ \delta[A \cup B] \leq K_A + K_B + \frac{2K_AK_B}{\eta}.\] Since $K_A, K_B \geq 1$ and $\eta \leq 1$, the result follows immediately.

(v) See \cite[Section 6.4]{tao-vu} for a proof and references to the original papers.

(vi) Suppose that $E(\lambda \cdot A, B) \geq \delta |A|^{3/2}|B|^{3/2}$. From the trivial estimates in (i) we see that $\delta^2 |A| \leq |B| \leq \delta^{-2}|A|$. By the Balog-Szemer\'edi-Gowers theorem there are sets $A' \subseteq A$ and $B' \subseteq B$ with $|A'| \geq c\delta^C |A|$, $|B'| \geq c\delta^C |B|$ such that $|\lambda \cdot A' + B'| \leq C\delta^{-C}|A'|$. By the Pl\"unnecke-Ruzsa inequality (cf. \cite[Cor 6.28]{tao-vu}) there is $A'' \subseteq A'$ such that
\[ |\lambda \cdot A'' + \lambda B'| \leq (C/\delta)^{CL}|A''|.\]
Since $\lambda \cdot(A'' + B') \subseteq \lambda \cdot A'' + \lambda B'$, this implies that
\[ |A'' + B'| \leq (C/\delta)^{CL}|A''|.\]
We clearly have \[ |A''| \geq (c\delta)^{CL}|B'| \geq (c\delta)^{CL}|A|.\] Thus from (iii) we 
obtain
\[ E(A,B) \geq E(A'',B') \geq (c\delta)^{CL}|A|^{3/2}|B|^{3/2}.\]
This implies the result. 
\endproof

\begin{proposition}[Structure theorem]\label{structure-theorem}
Let $A$ be an additive set and let $\epsilon,\epsilon' \in (0,1/2)$ be parameters. Let $L \in \N$ be fixed. Then there is a decomposition of $A$ as a disjoint union $A_1 \cup \dots \cup A_k \cup A_0$ such that
\begin{enumerate}
\item \textup{(Components are large)} $|A_i| \geq |A|/F_1(L,\epsilon)$ for $i=1,\dots,k$;
\item \textup{(Components are structured)} $\delta[A_i] \leq F_2(L,\epsilon,\epsilon')$ for $i = 1,\dots,k$;
\item \textup{(Different components do not communicate)} $E(\lambda_i \cdot A_i,\lambda_j \cdot A_j) \leq \epsilon' |A_i|^{3/2}|A_j|^{3/2}$ whenever $1 \leq i < j \leq k$ and whenever $\lambda_i, \lambda_j \in \{1,\dots,L\}$;
\item \textup{(Noise term)} $E(\lambda_0 \cdot A_0,\lambda \cdot A) \leq \epsilon |A|^3$ whenever $\lambda_0,\lambda \in \{1,\dots,L\}$.
\end{enumerate}
\end{proposition}
\emph{Remarks.} Property (i) guarantees that $k \leq F_1(L,\epsilon)$; that is, the ``complexity'' $k$ of the decomposition is bounded. We may take  
\[ F_1(L,\epsilon) = (C/\epsilon)^{CL^2}\] and
\[ F_2(L,\epsilon,\epsilon') = (C/\epsilon\epsilon')^{(C/\epsilon)^{CL^2}}\]
though the precise form of these bounds is not important for our application.

\proof Take $\eta := (c\epsilon)^{CL^2}$ and $\eta' := (c\epsilon')^{CL^2}$. If $C,c$ are chosen appropriately it will be enough to establish the proposition with (iii) replaced by 
\begin{equation} \tag*{$\mbox{(iii)}^{\prime}$}
E(A_i,A_j) \leq \eta' |A_i|^{3/2}|A_j|^{3/2}\end{equation} and (iv) replaced by \begin{equation} \tag*{$\mbox{(iv)}^{\prime}$} E(A_0,A) \leq \eta |A|^3.\end{equation} Statements (iii) and (iv) then follow automatically in view of Lemma \ref{basic-add-energy} (vi).

In this proof the reader should be particularly aware of the fact that the absolute constant $C$ may change from line to line. We begin by applying the Balog-Szemer\'edi-Gowers theorem iteratively. We will define a sequence of disjoint sets $B_1,B_2,\dots$. These having been defined, set $S_i := A \setminus (B_1 \cup \dots \cup B_i)$ (with the convention that $S_0 = A$). If, for some $i$, we have $E(A,S_i) \leq \eta |A|^3$ then we stop the iteration and set $A_0 := S_i$. If not then the Balog-Szemer\'edi-Gowers theorem informs us that there are sets $A' \subseteq A, S'_i \subseteq S_i$ with $|A'| \geq c\eta^{C}|A|$ and $|S'_i| \geq c\eta^C|S_i|$ such that $|A' -  S'_i| \leq C\eta^{-C} |A'|^{1/2}|S'_i|^{1/2}$. Set $B_{i+1} := S'_i$. Then by the Ruzsa triangle inequality we have $\delta[B_{i+1}] \leq C\eta^{-C}$. By Lemma \ref{basic-add-energy} (i) we have $|B_{i+1}| \geq c\eta^C|A|$. It follows that the iteration must stop after at most $s \leq C\eta^{-C}$ steps. 

Now the sets $B_i$ satisfy (i), (ii) and (iv)'. However (iii)' may fail, that is to say there may be additive communication between the sets $B_i$. If we do have $E(B_i,B_j) \geq \eta'|B_i|^{3/2}|B_j|^{3/2}$ for some $i \neq j$ then we simply replace $B_i$ and $B_j$ by the single set $B_i \cup B_j$, noting from Lemma \ref{basic-add-energy} (iv) that we have
\[ \delta[B_i \cup B_j] \leq 4\delta[B_i]\delta[B_j]/\eta'.\]
We then repeat if necessary.
It is clear that this process of ``agglomeration'' lasts no more than $s$ steps, in which time the $\delta[ \;]$ constants of all sets are still bounded by $(C/\eta\eta')^{C/\eta^C}$. This concludes the proof.\endproof

\emph{Remarks.} One can envisage various refined versions of this result, but we do not describe them in detail here for want of applications. Similar refinements were also discussed by Elekes and Ruzsa. Perhaps the most obvious step is to apply Freiman's theorem to each of the $A_i$, thereby placing $A \setminus A_0$ inside a union of multidimensional progressions $P_i$. One could easily ensure, by an agglomeration process similar to that used in the proof of Proposition \ref{structure-theorem}, that $E(P_i,P_j)$ is small when $i \neq j$. One might even go further, subdividing each $P_i$ into structured pieces (such as Bohr sets) such that $A$ looks pseudorandom on most of these pieces. By analogy with a result of the first author and Tao (\cite[Proposition 3.9]{green-montreal}), T.~Tao has suggested that such a result might be called a type of ``arithmetic regularity lemma''. Such a result would only be of use for qualitative applications -- such as that in the present paper -- as it would come with bounds of tower type.

\section{Structure, rectification and 3APs}

In this section we combine the structure theorem with a result of Bilu, Lev and Ruzsa \cite{blr}. This will first allow us, in Lemma \ref{lem4.3}, to place an upper bound on the number of 3APs in a set which has been decomposed as in Proposition \ref{structure-theorem}. We will then use that lemma to obtain an approximate structural result for subsets of $\Z/N\Z$ with close to the maximal number of 3APs. In the next section we will bootstrap that approximate result to an exact result.

The result of Bilu, Lev and Ruzsa to which we refer is a rectification lemma of the following type. The bounds stated below are those given in \cite{green-ruzsa}, which has the advantage of not requiring Freiman's theorem for its proof.

\begin{theorem}[Rectification lemma]\label{sumset-rectify}
Suppose that $N$ is a prime. Let $B \subset \Z/N\Z$ be a set with $|B|=\beta N$ such that $\delta[B] \leq K$. Suppose that $\beta \leq (16K)^{-12K^2}$. Then there is $d \in (\Z/N\Z)^*$ such that $d\cdot B$ is contained in an interval of length at most 
\[ 12\beta^{1/4K^2}\sqrt{\log(1/\beta)} N. \]
\end{theorem}

The next lemma provides a bound for the number of 3APs in $A_1 \times A_2 \times A_3$ in terms of additive energies. 
\begin{lemma}[Bounding 3APs using the additive energy]\label{lem3.1}
Suppose that $A_1,A_2,A_3$ are three subsets of an abelian group. Then 
\[ T_3(A_1,A_2,A_3)^6 \leq |A_1||A_2||A_3|E(2\cdot A_2, A_3)E(A_1,A_3)E(A_1,2\cdot A_2).\] 
\end{lemma}
\proof For each $y \in A_2$, let $m(y)$ denote the number of pairs $(x,z) \in A_1 \times A_3$ such that $(x+z)/2 = y$. Thus $\sum_y m(y) = T_3(A_1,A_2,A_3)$. Now $\sum_y m(y)^2$ is at most the number of solutions to $x + z = x' + z'$ with $x,x' \in A_1$ and $z,z' \in A_3$, which is precisely $E(A_1,A_3)$. Thus by the Cauchy-Schwarz inequality we have
\[ E(A_1,A_3) = \sum_y m(y)^2 \geq \frac{1}{|A_2|} \big(\sum_y m(y)\big)^2 = T_3(A_1,A_2,A_3)^2/|A_2|.\]
There are two similar lower bounds for $E(2\cdot A_2,A_3)$ and $E(A_1,2\cdot A_2)$, which may be proved in exactly the same way. Multiplying the three bounds together gives the result.
\endproof

\begin{lemma}[Bounding 3APs in sets]\label{lem4.3}
Suppose that $N$ is a prime and that $\delta \in (0,1)$. Then there are $\epsilon,\epsilon' \leq C\delta^{-C}$ and a constant $c_{\delta} > 0$ with the following property. Let $A \subseteq \Z/N\Z$ be any set of cardinality $n$ satisfying $c_{\delta}^{-1} \leq n \leq c_{\delta}N$. Apply Proposition \ref{structure-theorem} with parameters $\epsilon,\epsilon'$ to obtain a decomposition $A =  A_1 \cup \dots \cup A_k \cup A_0$ satisfying conditions (i), (ii), (iii) and (iv) of that proposition with $L = 2$. Then, writing $n_i := |A_i|$ for $i = 1,\dots,k$, we have
\[ T_3(A) \leq \frac{1}{2}\sum_{i=1}^k n_i^2 + \delta n^2.\]
\end{lemma}
\proof Take $\epsilon := (\delta/9)^3$ and define $\epsilon' := (\delta/3F_1(2,\epsilon))^6$, where $F_1$ is the function occurring in Proposition \ref{structure-theorem}. Recall that, in particular, $F_1$ provides a bound for the ``complexity'' $k$ of the decomposition $A = A_1 \cup \dots \cup A_k \cup A_0$.

We of course have
\[ T_3(A) = \sum_{0 \leq i_1,i_2,i_3 \leq k} T_3(1_{A_{i_1}},1_{A_{i_2}},1_{A_{i_3}}),\]
a sum which we split into three parts $S_1$, $S_2$ and $S_3$. $S_1$ is the contribution from the terms $i_1 = i_2 = i_3 > 0$, $S_2$ is the contribution from the terms where some $i$ equals zero, and $S_3$ is the contribution from the remaining terms, those with $i_1,i_2,i_3 > 0$ and not all equal.

\emph{Bounding $S_1$.} Clearly
\[ S_1 = \sum_{i=1}^k T_3(A_i).\]
Now the set $A_i$ satisfies $\delta[A_i] \leq F_2(2,\epsilon,\epsilon')$, and so (provided that $c_{\delta}$ is chosen sufficiently small) Theorem \ref{sumset-rectify} guarantees that some dilate $A'_i$ of $A_i$ is contained in a translate of $[0,\lfloor N/2\rfloor] \subseteq \Z/N\Z$. We may associate to this set the corresponding set $A^*_i \subseteq [0,\lfloor N/2\rfloor]$ of integers, and it is clear that
\[ T_3(A_i) = T_3(A'_i) = T_3(A^*_i).\]
It therefore follows from Theorem \ref{max3APInt} that
\[ T_3(A_i) \leq \lceil n_i^2/2\rceil \] and hence that
\[ S_1 \leq \frac{1}{2}\sum_{i=1}^k n_i^2 + \delta n^2/3 .\] (Note that the \emph{lower} bound $n \geq c_{\delta}^{-1}$ is required here.)

\emph{Bounding $S_2$.} We have
\[ S_2 \leq T_3(1_{A_0}, 1_A, 1_A) + T_3(1_A, 1_{A_0}, 1_A) + T_3(1_A, 1_A, 1_{A_0}).\]
By Lemma \ref{lem3.1} and the property of Proposition \ref{structure-theorem} (iv) each term is bounded by $\epsilon^{1/3}n^2$, and so $S_2 \leq \delta n^2/3$.

\emph{Bounding $S_3$.} Provided that $i_1,i_2,i_3 > 0$ and are not all equal, Lemma \ref{lem3.1} tells us that 
\[ T_3(1_{A_{i_1}}, 1_{A_{i_2}}, 1_{A_{i_3}}) \leq \epsilon^{\prime 1/6} n_{i_1}^{2/3} n_{i_2}^{2/3} n_{i_3}^{2/3}.\]
Summing over $i_1,i_2,i_3$ and using the fact that $\sum_{i=1}^k n_i^{2/3} \leq k^{1/3}n^{2/3}$ (a consequence of H\"older's inequality) we have
\[ S_3 \leq \epsilon^{\prime 1/6} k n^2 \leq \delta n^2/3.\]
Putting together these three estimates for $S_1,S_2$ and $S_3$ leads to the result.\endproof

We now derive our approximate structure theorem for sets with close to the maximal number of 3APs.

\begin{lemma}\label{approx-struct}
There is an absolute constant $c > 0$ with the following property. Suppose that $n \leq cN$, and that $A \subseteq \Z/N\Z$ is a set with $|A| = n$ and $T_3(A) \geq 0.96 M_3(n,N)$ \textup{(}that is, $A$ has close to the maximal number of 3APs for subsets of $\Z/N\Z$ of size $n$\textup{)}. Then there is some dilate of $A$, at least 95\% of whose elements lie in a translate of the interval $[-N/24,N/24]$.
\end{lemma}
\proof If $|A| \leq C$ then there is a dilate and translate of $A$, \emph{all} of whose points lie in $[-N/24,N/24]$, by a standard application of the pigeonhole principle. If this is not the case then we may apply Lemma \ref{lem4.3} with $\delta = 1/200$. Provided $c$ is chosen sufficiently small this provides a decomposition of $A$ as $A_1 \cup \dots \cup A_k \cup A_0$ where $|A_i| = n_i$ and the number of 3APs in $A$, $T_3(A)$, is bounded by
\[ T_3(A) \leq \frac{1}{2}\sum_{i = 1}^k n_i^2 + \frac{1}{200}n^2.\]
Now, since $M_3(n,N) \geq n^2/2$, we have that $T_3(A) \geq 0.48n^2$, and therefore
\[ \sum_{i=1}^k n_i^2 \geq 19n^2/20.\]
Supposing without loss of generality that $n_1$ is the largest of the $n_i$ we see immediately that
\[ n_1 n \geq n_1 \sum_{i=1}^k n_i \geq 19n^2/20.\]
This implies that $n_1 \geq 19n/20$, which is to say that $A_1$ accounts for at least 95\% of the set $A$. To complete the proof of the lemma, then, we need only show that some dilate of $A_1$ lies in a translate of $[-N/24,N/24]$. This, however, is an immediate consequence of condition (ii) of Proposition \ref{structure-theorem} together with Theorem \ref{sumset-rectify}, provided that $c$ is chosen sufficiently small.\endproof

\section{From approximate structure to exact structure}

Our objective in this section is to conclude the proof of Theorem \ref{mainthm}. In view of Lemma \ref{approx-struct} it suffices to establish the following result.

\begin{lemma} 
Suppose that $A \subseteq \Z/N\Z$ has $|A| = n$ and that at least 95\% of the elements of $A$ lie in $[-N/24,N/24]$. Then $T_3(A) \leq \lceil n^2/2\rceil$ with equality if and only if $A$ is an affine copy of one of the sets $E(k,m),F(k,m)$, considered as a subset of $\Z/N\Z$. 
\end{lemma}
\proof Write \[ A_0 := A \cap [-N/24,N/24],\] \[ A_1 := A \cap ([-N/8,N/8] \cup [3N/8,5N/8]) \] and
\[ A_2 := A \setminus A_1.\]
Write $n_i := |A_i|$, $i = 1,2$.
One may check that any 3AP with at least \emph{two} points in $A_0$ must be entirely contained in $A_1$. Now $2 \cdot A_1 \subseteq [-N/4,N/4]$, and thus $A_1$ is Freiman-isomorphic to a set of integers and $T_3(A_1) \leq \lceil n_1^2/2\rceil$.

Any 3AP in $A$ that is not entirely contained in $A_1$ has at least one point in $A_2$, and can have at most one point in $A_0$. The number of such 3APs is therefore at most $6|A_2||A \setminus A_0| \leq \lfloor 3n_2 n/10\rfloor$.
Now if $n_2 \neq 0$ then we have
\[ \frac{n_1^2+1}{2} + \frac{3n_2 n}{10} < \frac{n^2}{2}\] and therefore
\[ T_3(A) < \lceil n^2/2\rceil.\]
If $n_2 = 0$ then $A = A_1$ and hence, as we have already seen, $2 \cdot A \subseteq [-N/4,N/4]$. Thus $T_3(2 \cdot A) \leq \lceil n^2/2\rceil$ with equality if and only if $2 \cdot A$ is an affine copy of one of the sets $E(k,m), F(k,m)$. The same is therefore true of $A$.\endproof

\section{Some remarks on Croot's functions $M_3(\alpha),m_3(\alpha)$.}
Recall 
that 
\[ M_3(\alpha) := \lim_{\substack{N \rightarrow \infty \\ N \mbox{\scriptsize \textup{prime}}}} \frac{M_3(\lfloor \alpha N\rfloor, N)}{N^2}\qquad \mbox{and}
\qquad m_3(\alpha) := \lim_{\substack{N \rightarrow \infty \\ N \mbox{\scriptsize \textup{prime}}}} \frac{m_3(\lfloor \alpha N\rfloor, N)}{N^2}.\]
Our main theorem together with the relation $m_3(\alpha) + M_3(1-\alpha) = 1 - 3\alpha + 3\alpha^2$ implies that
\[ M_3(\alpha) = \frac{\alpha^2}{2} \] for $\alpha<c$ and 
\[ m_3(\alpha) = \frac{1}{2} - 2\alpha + \frac{5}{2}\alpha^2 \]
for $\alpha>1-c$. We note here that we must certainly have $c \leq 1/3$ for these bounds to hold, for once $\alpha > 1/3$ one may choose particular values of $k$ and $m$ so that the set $E(k,m)$ exploits `wrap-around' in $\Z/N\Z$ to increase the 3AP-count quite significantly. In particular, the sets $E(k,m) \subset \Z/N\Z$ are in general not rectifiable for $\alpha>1/3$.

More specifically, choosing $k \approx (3n-N)/6$ one can (using the complement $E(k,m)^c$) obtain the bound 
\begin{align*} m_3(\alpha) \leq (2 - 12\alpha + 21\alpha^2)/12 \end{align*}
for $1/3 \leq \alpha \leq 2/3$. For $\alpha=1/2$, this leads to the bound
\begin{align} m_3(1/2) \leq 5/48. \label{Uhalfbound} \end{align}
By contrast, arithmetic progressions and random sets of density $1/2$ lead only to the trivial $m_3(1/2) \leq 1/8$. The bound \eqref{Uhalfbound} seems to be the best we currently have for Croot's function $m_3$ at $\alpha=1/2$, though there is perhaps insufficient evidence to risk conjecturing that it represents the true state of affairs.

In another paper \cite{croot2}, Croot makes the pleasant observation that
\[ m_3(\alpha \beta) \leq m_3(\alpha)m_3(\beta).\]
To prove this he takes sets $A,B$ with $|A| \approx \alpha N$, $T_3(A) \approx m_3(\alpha)N^2$, $|B| \approx \beta N$, $T_3(B) \approx m_3(\beta) N^2$ and looks at the intersection $A \cap (\lambda B + \mu)$ for randomly chosen $\lambda \in (\Z/N\Z)^*$ and $\mu \in \Z/N\Z$. Using the first and second moment methods he shows that with positive probability one has $|A \cap (\lambda B + \mu)| \approx \alpha \beta N$ and $T_3(A \cap (\lambda B + \mu)) \leq m_3(\alpha)m_3(\beta)N^2+o(N^2)$, thereby establishing the result. (One can show in almost identical fashion that $M_3(\alpha\beta) \geq M_3(\alpha)M_3(\beta)$.)

Using this observation one may find a cutoff density below which one can be certain that the upper bound for $m_3(\alpha)$ given by the sets $E(k,m)$ (or their complements, rather) is no longer sharp. This is certainly the case for $\alpha < 2(7+2\sqrt{6})/75 \approx 0.3173$. In other words, once we are below a certain density, a ``2-dimensional'' set consisting of the intersection of a set $E(k,m)^c$ with an affine image of itself does rather better than any single set $E(k,m)^c$. 

The appearance of such multidimensional examples is perhaps not surprising in view of the fact that the best known constructions of sets with very few 3APs (with $\alpha$ very small) come from variants of the Behrend example \cite{behrend}, which is a sort of projection of the set of lattice points on a high-dimensional sphere.

\label{miscell-sec}

\section{Acknowledgement}

We would like to thank Tom Sanders for helpful conversations, especially concerning Proposition \ref{structure-theorem}.

     \end{document}